\input amstex
\documentstyle{amsppt}
\magnification = 1200
\pageheight{7.5in}
\vcorrection{-.12in}
\NoBlackBoxes
\NoRunningHeads
\TagsOnRight
\topmatter
\title The Maximum Principle for the Bergman Space \\and the
M\"obius Pseudodistance for the Annulus
\endtitle
\rightheadtext{The Bergman Space Maximum Principle}
\author Alexander Schuster
\endauthor
\address Department of Mathematics, San Francisco State University,
San Francisco, CA 94132.
\endaddress

\affil San Francisco State University
\endaffil
\email schuster\@sfsu.edu
\endemail
\abstract It is shown that the formula for the M\"obius
pseudodistance for the annulus yields better estimates than
previously known for the constant in the Bergman space maximum
principle.
\endabstract
\endtopmatter
\document
{\bf \S 1. Introduction.}

The Bergman space $A^2$ is the set of functions analytic in the unit
disk $\Bbb D=\{z\in\Bbb C: |z|<1\}$ with
$$\|f\|=\left\{\frac{1}{\pi} \int_{\Bbb D}
|f(z)|^2dA(z)\right\}^{\frac{1}{2}}<\infty,$$ where $dA$ denotes
Lebesgue area measure. An important result in the theory of Bergman
spaces is the so-called Korenblum maximum principle, which we shall
also refer to as the Bergman space maximum principle.

\proclaim{Theorem 1} There exists $c\in (0,1)$ with the property
that whenever $f$ and $g$ are functions that are analytic in $\Bbb
D$ and satisfy $|f(z)|\leq |g(z)|$ for $c<|z|<1$, then
$\|f\|\leq\|g\|$.
\endproclaim

First conjectured by Korenblum [5], the maximum principle was proved
by Hayman [2] with $c=0.04$. Hinkkanen [3] later
 improved upon Hayman's constant by showing that the result holds for $c=0.157...$.
 Moreover, he proved that it is valid more generally in
 the $L^p$ Bergman
space, where  $1\leq p<\infty$. The value of the best constant, even
for $p=2$, remains a mystery. Setting $f(z)=c$ and $g(z)=z$ shows
that $c<1/\sqrt{2}=0.707...$. Wang [6] has recently shown that
$c<0.69472...$ by considering a certain class of singular inner
functions. Our contribution in this paper is to show that the
Korenblum maximum principle holds with $c=0.21$.

The paper is organized as follows. In section 2, much of which
follows Hinkkanen [3], we introduce the notation and obtain some
preliminary estimates. In section 3 we discuss the M\"obius
pseudodistance of the annulus and apply it to obtain the necessary
estimates for our proof of the maximum principle. Section 4 consists
of the application of our ideas to the setting of
 the Fock space of entire functions.
$$ $$

{\bf \S 2. Preliminaries.}

Let $c\in(0,1)$ be a constant to be determined later, and suppose
that $f$ and $g$ are functions analytic in $\Bbb D$ satisfying
$|f(z)|\leq |g(z)|$ for $z\in A(c,1)$, where $A(r_1,r_2)=\{z\in\Bbb
C:r_1<|z|<r_2\}$. Our goal is to show that
$$\int_{D(c)}(|f(z)|^2-|g(z)|^2)dA(z)\leq \int_{A(c,1)}(|g(z)|^2-|f(z)|^2)dA(z),\tag1$$
where $D(c)$ is the open disk of radius $c$ centered at the origin.

Define the function $\omega=f/g$, which by hypothesis is analytic and satisfies $|\omega|\leq1$ in  $A(c,1)$. In fact, we may assume,
without loss of generality, that $|\omega|<1$ in $A(c,1)$, since otherwise $|f|=|g|$ in $\Bbb D$ and the result
holds trivially. We may likewise assume that $f$ is not identically equal to $0$.

For $c<\rho<1$, choose $\zeta_\rho$ such that $|\zeta_\rho|=\rho$
and $|\omega(\zeta_\rho)|=\sup\{|\omega(z)|: |z|=\rho\}$, and define
$\omega_\rho=\omega(\zeta_\rho)$. We will assume, without loss of
generality, that $\zeta_\rho=\rho$.

Let $0<r<c<\rho<1$. Following Hinkkanen, we use the inequality
$|\alpha|^2-|\beta|^2\leq 2|\alpha^2-\alpha\beta|$ and the fact
that, for a subharmonic function $h$ and $0<r_1<r_2<1$,
$\int_0^{2\pi}h(r_1e^{i\theta})d\theta\leq\int_0^{2\pi}h(r_2e^{i\theta})d\theta$,
to obtain
$$\align
\int_0^{2\pi}(|f(re^{i\theta})|^2&-|g(re^{i\theta})|^2)d\theta \leq
\int_0^{2\pi}(|f(re^{i\theta})|^2-|\omega_\rho g(re^{i\theta})|^2)
d\theta\\
&\leq2\int_0^{2\pi}|f^2(re^{i\theta})-
\omega_\rho f(re^{i\theta})g(re^{i\theta})|d\theta\\
&\leq2\int_0^{2\pi}|f^2(\rho e^{i\theta})-\omega_\rho
f(\rho e^{i\theta})g(\rho e^{i\theta})|d\theta\\
&=2\int_0^{2\pi}|\omega(\rho e^{i\theta})|\frac{|\omega(\rho e^{i\theta})-\omega_\rho|}{1-|\omega(\rho e^{i\theta})|^2}(|g(\rho e^{i\theta})|^2-|f(\rho e^{i\theta})|^2)d\theta\\
&\leq2\gamma(\rho)\int_0^{2\pi}(|g(\rho e^{i\theta})|^2-|f(\rho e^{i\theta})|^2)d\theta,
\endalign$$
where
$$\gamma(\rho)
=\sup\left\{\frac{|\omega(z)-\omega_\rho|}{1-|\omega(z)|^2}:
|z|=\rho\right\}.$$ Multiplying both sides of this inequality by $r$
and then integrating from $0$ to $c$ with respect to $r$ yields
$$\int_{D(c)}(|f(z)|^2-|g(z)|^2)dA(z)\leq c^2\gamma(\rho)\int_0^{2\pi}
(|g(\rho e^{i\theta})|^2-|f(\rho e^{i\theta})|^2)d\theta.$$ Multiply
both sides of this last inequality by $\rho(\gamma(\rho))^{-1}$ and
integrate from $c$ to $1$ with respect to $\rho$ to arrive at
$$\align
\int_{D(c)}(|f(z)|^2-|g(z)|^2)&dA(z)\\&\leq
c^2\left(\int_c^1\gamma(\rho)^{-1}\rho d\rho\right)^{-1}
\int_{A(c,1)}(|g(z)|^2-|f(z)|^2)dA(z).\tag1
\endalign$$

It remains to find $c$ such that the quantity in front of the last
integral is bounded above by $1$. To achieve this goal we will need
to find a suitable bound for $\gamma(\rho)$, which we proceed to do
in the next section.
$$ $$

{\bf \S 3. The M\"obius pseudodistance }

To estimate $\gamma(\rho)$ we recall that the pseudohyperbolic
distance $d$ between two points $\alpha,\beta\in\Bbb D$ is given by
the formula
$$d(\alpha,\beta)=\left|\frac{\alpha-\beta}{1-\overline{\alpha}\beta}\right|.$$
By the identity
$$\frac{|\alpha-\beta|}{1-|\alpha|^2}=\frac{d(\alpha,\beta)}
{\sqrt{1-d^2(\alpha,\beta)}}\frac{\sqrt{1-|\beta|^2}}{\sqrt{1-|\alpha|^2}},$$
we see that for $|z|=\rho$,
$$\frac{|\omega(z)-\omega_\rho|}{1-|\omega(z)|^2}=\frac{d(\omega(z),\omega_\rho)}
{\sqrt{1-d^2(\omega(z),\omega_\rho)}}\frac{\sqrt{1-|\omega_\rho|^2}}
{\sqrt{1-|\omega(z)|^2}}\leq \frac{d(\omega(z),\omega_\rho)}
{\sqrt{1-d^2(\omega(z),\omega_\rho)}}.$$
It will therefore behoove us to find an estimate for
$d(\omega(z),\omega_\rho)$. To do this, we consider the following
notion of the distance between two points.

 The M\"obius pseudodistance for a domain $D\subset\Bbb C$ is
defined by the equation
$$c_D^*(a,z)=\sup\{d(\omega(a),\omega(z)):\omega\in Hol(D,\Bbb D)\},$$ where
$a,z\in D$ and $Hol(A,B)$ denotes the set of analytic functions from
$A$ to $B$. (See the book of Jarnicki and Pflug [4] for an excellent
survey of this topic.) A basic property of the M\"obius
pseudodistance is its invariance with respect to biholomorphic maps.
Namely, if $\psi:D\to\Omega$ is a biholomorphism, then
$$c_\Omega^*(\psi(a),\psi(z))=c_D^*(a,z).$$ Since the pseudohyperbolic
metric is M\"obius invariant,
$$c_D^*(a,z)=\sup\{|\omega(z)|:\omega\in Hol(D,\Bbb D),
\omega(a)=0\}.$$

It is shown in [1] that for the annulus $P=A(\frac{1}{R},R)$ and
$\frac{1}{R}<a<R$,
$$c^*_{P}(a,z)=\frac{f(\frac{1}{a},-|z|)}{R|z|}|f(a,z)|,$$
where $$f(a,z)=\left(1-\frac{z}{a}\right)\prod_{n=1}^\infty\frac{
(1-\frac{z}{a}R^{-4n})(1-\frac{a}{z}R^{-4n})}{
(1-azR^{2-4n})(1-\frac{1}{az}R^{2-4n})}.$$

 From the invariance of the
M\"obius pseudodistance under biholomorphic maps we obtain the
formula
$$\align &c_{A(c,1)}^*(\rho,z)=c_{A(\sqrt{c},\frac{1}{\sqrt{c}})}^*
(\frac{\rho}{\sqrt{c}},\frac{z}{\sqrt{c}})=
\frac{c\,f(\frac{\sqrt{c}}{\rho},-\frac{|z|}{\sqrt{c}})
|f(\frac{\rho}{\sqrt{c}},\frac{z}{\sqrt{c}})|}{|z|}\\
&=\frac{c}{|z|}\left(1+\frac{\rho|z|}{c}\right)\left(1-\frac{z}{\rho}\right)\left|\prod_{n=1}^\infty
\frac{ (1+\rho|z|c^{2n-1})(1+\frac{1}{\rho|z|}c^{2n+1})
(1-\frac{z}{\rho}c^{2n})(1-\frac{\rho}{z}c^{2n})}
{(1+\frac{|z|}{\rho}c^{2n-1})(1+\frac{\rho}{|z|}c^{2n-1}) (1-\rho
zc^{2n-2})(1-\frac{1}{\rho z}c^{2n})}\right|.
\endalign$$
In particular,
$$c_{A(c,1)}^*(\rho,\rho e^{i\theta})=\frac{c}{\rho}\left(1+\frac{\rho^2}{c}\right)\sqrt{2(1-\cos\theta)}
\prod_{n=1}^\infty f_n(\rho,c)g_n(\rho,c,\theta),$$ where
$$f_n(\rho,c)=\frac{(1+\rho^2c^{2n-1})(1+\rho^{-2}c^{2n+1})}{(1+c^{2n-1})^2},$$
and
$$g_n(\rho,c,\theta)=\frac{1-2c^{2n}\cos\theta+c^{4n}}
{\sqrt{1-2\rho^2c^{2n-2}\cos\theta+\rho^4c^{4n-4}}
\sqrt{1-2\rho^{-2}c^{2n}\cos\theta+\rho^{-4}c^{4n}}}.$$ Note that
$$f_n(\rho,c)\leq \frac{1+c^{2n+1}}{1+c^{2n-1}}.$$
A calculus argument shows that $$g_n(\rho,c,\theta)\leq
g_n(\rho,c,0)=\frac{(1-c^{2n})^2}{(1-\rho^2c^{2n-2})(1-\rho^{-2}c^{2n})}
\leq\frac{1-c^{2n}}{1-c^{2n-2}},$$ and so
$$f_n(\rho,c)g_{n+1}(\rho,c)\leq\frac{(1+c^{2n+1})(1-c^{2n+2})}{(1+c^{2n-1})(1-c^{2n})},$$
which is bounded above by $1$ for every $n$ if $c<1/4$. Therefore
$$\prod_{n=6}^\infty f_n(\rho,c)g_n(\rho,c)\leq
\left(\frac{1-c^{12}}{1-c^{10}}\right)\prod_{n=6}^\infty
f_n(\rho,c)g_{n+1}(\rho,c) \leq \frac{1-c^{12}}{1-c^{10}}.$$

On the other hand, a long and extremely tedious calculation shows
that
$$\sqrt{2(1-\cos\theta)}\prod_{n=1}^5 g_n(\rho,c,\theta)\leq
2\prod_{n=1}^5g_n(\rho,c,\pi)=2\prod_{n=1}^5
\frac{(1+c^{2n})^2}{(1+\rho^2c^{2n-2})(1+\rho^{-2}c^{2n})}.$$

 Putting everything together, we have
 $$c_{A(c,1)}^*(\rho,\rho e^{i\theta})\leq F(\rho,c),$$
 where
 $$F(\rho,c)=2\frac{c}{\rho}\left(1+\frac{\rho^2}{c}\right)
 \left(\frac{1-c^{12}}{1-c^{10}}\right)\prod_{n=1}^5
 \frac{(1+\rho^2c^{2n-1})(1+\rho^{-2}c^{2n+1})(1+c^{2n})^2}{(1+c^{2n-1})^2
 (1+\rho^2c^{2n-2})(1+\rho^{-2}c^{2n})}.$$
Then
$$\sup_{|z|=\rho}c^*_{A(c,1)}(\rho,z)\leq F(\rho,c),$$ and so
$$\gamma(\rho)\leq \sup_{|z|=\rho}\frac{c^*_{A(c,1)}(\rho,z)}
{\sqrt{1-(c^*_{A(c,1)}(\rho,z))^2}}\leq\frac{F(\rho,c)}{\sqrt{1-F^2(\rho,c)}}.$$
We combine this with (1) to obtain
$$\align
\int_{D(c)}(|f(z)|^2&-|g(z)|^2)dA(z)\\&\leq
c^2\left(\int_c^1\frac{\rho\sqrt{1-F^2(\rho,c)}}{F(\rho,c)}
d\rho\right)^{-1} \int_{A(c,1)}(|g(z)|^2-|f(z)|^2)dA(z).
\endalign$$

A calculation involving Mathematica now shows that when $c=0.21$,
the  quantity in front of the last integral is less than $1$, and so
the proof of the result is complete.
$$ $$

{\bf \S 4. The Fock space}

The Fock space $F$ is the set of entire functions with
$$\|f\|_F=\left\{\int_{\Bbb C}|f(z)|^2e^{-|z|^2}dA(z)\right\}
^{\frac{1}{2}}<\infty.$$

We prove the following analogue of the Korenblum maximum principle
in this setting.

\proclaim{Theorem 2} There is a positive constant $c$ with the
property that whenever $f$ and $g$ are entire functions satisfying
$|f(z)|\leq |g(z)|$ for $|z|> c$, then $\|f\|_F\leq\|g\|_F$.
\endproclaim

Note that the proof of the Bergman space maximum principle can be
modified to give the result that there is a constant $c$ with the
property that whenever $f$ and $g$ are analytic in $\Bbb D$ with
$|f(z)|\leq |g(z)|$ for $c<|z|<1$, then
$$\int_{\Bbb D}|f(z)|^2e^{-|z|^2}dA(z)\leq\int_{\Bbb
D}|g(z)|^2e^{-|z|^2}dA(z).$$ In particular, the maximum principle
for the Fock space is a consequence of the maximum principle for the
Bergman space. However, the proof of the former is significantly
easier than the proof of the latter, and in addition, we obtain a
better constant ($c=0.54$), so we give the proof below.

All of the notation will be exactly the same as in our proof of
Theorem 1. The first part of the argument will be suitably modified
to obtain the inequality
$$\align
\int_{D(c)}&(|f(z)|^2-|g(z)|^2)e^{-|z|^2}dA(z) \\&\leq
(1-e^{-c^2})\left(\int_c^\infty\gamma(\rho)^{-1} \rho
d\rho\right)^{-1}\int_{A(c,\infty)}
(|g(z)|^2-|f(z)|^2)e^{-|z|^2}dA(z). \endalign$$
Again we have
$$\gamma(\rho)\leq \sup_{|z|=\rho}\frac{c^*_{A(c,\infty)}(z,\rho)}
{\sqrt{1-(c^*_{A(c,\infty)}(z,\rho))^2}},$$ but here this quantity
is easier to estimate. To this end, we define the function $H:\Bbb
D^*\to\Bbb D$, where $\Bbb D^*$ is the punctured unit disk, by the
equation $H(\eta)=\omega(c/\eta)$. This is analytic in $\Bbb D^*$,
and since $\omega$ is bounded in $A(c,\infty)$, the function $H$ can
be extended to be analytic on all of $\Bbb D$. By the Schwarz-Pick
lemma, we then have for $|z|=\rho$,
$$c^*_{A(c,\infty)}(z,\rho)
\leq d\left(\frac{c}{z},\frac{c}{\zeta_\rho}\right)
\leq\frac{2c\rho}{\rho^2+c^2},$$ and so
$$\gamma(\rho)\leq \frac{2c\rho}{\rho^2-c^2}.$$ Therefore,
$$\align&\int_{D(c)}(|f(z)|^2-|g(z)|^2)e^{-|z|^2}dA(z)\\&
\leq 2c(1-e^{-c^2})\left(\int_c^\infty e^{-\rho^2}
(\rho^2-c^2)d\rho\right)^{-1}\int_{A(c,\infty)}
(|g(z)|^2-|f(z)|^2)e^{-|z|^2}dA(z).\endalign$$

A calculation involving Mathematica shows that when $c=0.54$, the
quantity in front of the last integral is less than $1$.

\demo{Acknowledgement} The author is grateful to Dror Varolin for
helpful comments and observations.
\enddemo

\enddocument